\documentclass[11pt,twoside, reqno]{amsart}
\date{26 April 2018}
\usepackage{latexsym,amsmath,amsthm,amsfonts,amscd,amssymb}
\usepackage[utf8]{inputenc}
\usepackage[all]{xy}
\usepackage{graphics}
\usepackage{lscape}
\usepackage{array}
\usepackage{framed}
\usepackage{mathtools}
\usepackage{hyperref}\hypersetup{colorlinks}

\usepackage{color} 

\definecolor{darkred}{rgb}{1,0,0} 
\definecolor{darkgreen}{rgb}{0,0.8,0}
\definecolor{darkblue}{rgb}{0,0,1}

\hypersetup{colorlinks,
linkcolor=darkblue,
filecolor=darkgreen,
urlcolor=darkred,
citecolor=darkgreen}

\theoremstyle{plain}  
\newtheorem{theorem}{Theorem}[section]

\newtheorem*{theorem*}{Theorem}

\newtheorem{corollary}[theorem]{Corollary}

\newtheorem{proposition}[theorem]{Proposition}
\newtheorem{conjecture}[theorem]{Conjecture}
\newtheorem{tech-lemma}[theorem]{Technical Lemma}
\newtheorem{definition}[theorem]{Definition}
\theoremstyle{remark}

\newtheorem{remark}[theorem]{Remark}
\newtheorem*{remark*}{Remark}

\newtheorem*{claim*}{Claim}
\numberwithin{equation}{section}
\renewcommand{\leq}{\leqslant}
\renewcommand{\le}{\leqslant}
\renewcommand{\geq}{\geqslant}

\newcommand{\Z}{\mathbb{Z}}
\newcommand{\CC}{\mathbb{C}}

\newcommand{\cM}{{\mathcal M}}

\newcommand{\rmin}{\mathrm{top}}

\newcommand{\U}{\mathrm{U}}

\newcommand{\G}{\mathrm{G}}

\newcommand{\SO}{\mathrm{SO}}
\newcommand{\Or}{\mathrm{O}}
\newcommand{\Sp}{\mathrm{Sp}}

\newcommand{\Sym}{\mathrm{Sym}}

\newcommand{\mtrx}[1]{\left (\begin{matrix}#1\end{matrix}\right)}
\newcommand{\smtrx}[1]{\left (\begin{smallmatrix}#1\end{smallmatrix}\right)}

\DeclareMathOperator{\Hom}{Hom}

\DeclareMathOperator{\Rep}{\mathcal{R}}

\newcommand{\R}{\mathbb{R}}

\newcommand{\C}{\mathbb{C}}
\newcommand{\noi}{\noindent}

\let\oldmarginpar\marginpar
\renewcommand\marginpar[1]{\oldmarginpar{\tiny\bf\begin{flushleft} #1
\end{flushleft}}}

\begin{document}

%
%

\title[Exotic components of $\SO(p,q)$ surface group representations]{Exotic components of $\SO(p,q)$ surface group representations, and their Higgs bundle avatars}
%
%


\author[M. Aparicio]{Marta Aparicio-Arroyo}
\address{Raet | HR software and services\newline\indent
Avenida de Bruselas 7, 28108 Alcobendas, Madrid, Spain}
\email{Marta.Aparicio@raet.com}

\author[S. Bradlow]{Steven Bradlow}
\address{Department of Mathematics, University of Illinois at Urbana-Champaign\newline\indent
Urbana, IL 61801, USA}
\email{bradlow@math.uiuc.edu}

\author[B. Collier]{Brian Collier}
\address{Department of Mathematics, University of Maryland\newline\indent College Park, MD 20742, USA}
\email{briancollier01@gmail.com}

\author[O. García-Prada]{Oscar García-Prada}
\address{Instituto de Ciencias Matem\'aticas,  CSIC-UAM-UC3M-UCM, \newline\indent 
Nicol\'as Cabrera, 13--15, 28049 Madrid, Spain}
\email{oscar.garcia-prada@icmat.es}

\author[P.~B.\ Gothen]{Peter B.\ Gothen}
\address{Centro de Matemática da Universidade do Porto, 
	\newline\indent Faculdade de Ci\^encias da Universidade do Porto, 
    \newline\indent Rua do Campo Alegre s/n, 4169-007 Porto, Portugal}
\email{pbgothen@fc.up.pt}

\author[A. Oliveira]{André Oliveira}
\address{Centro de Matemática da Universidade do Porto,
     \newline\indent Faculdade de Ci\^encias da Universidade do Porto, 
     \newline\indent
     Rua do Campo Alegre s/n, 4169-007 Porto, Portugal \newline\indent \textit{On leave from:}\newline\indent Departamento de Matemática, Universidade de Trás-os-Montes e Alto Douro, UTAD,\newline\indent
Quinta dos Prados, 5000-911 Vila Real, Portugal}
\email{andre.oliveira@fc.up.pt\newline\indent agoliv@utad.pt}


\begin{abstract} For semisimple Lie groups, moduli spaces of Higgs bundles on a Riemann surface correspond to representation varieties for the surface fundamental group.  In many cases, natural topological invariants label connected components of the moduli spaces.  Hitchin representations into split real forms, and maximal representations into Hermitian Lie groups, are the only previously know cases where natural invariants do not fully distinguish connected components. In this note we announce the existence of new such exotic components in the moduli spaces for the groups $\SO(p,q)$ with $2<p< q$. These groups lie outside formerly know classes of groups associated with exotic components.

\bigskip
\noindent \scalebox{.8}{R\'ESUM\'E.} Pour les groups de Lie semisimples, les espaces de modules de fibrés de Higgs sur une surface de Riemann
sont en correspondance avec les variétés de représentations du groupe fondamental de la surface. Pour beaucoup de groupes,  les invariants topologiques naturels  distinguent les composantes connexes de l'espace de modules.  Les représentations de Hitchin dans un groupe réel déployé et
des représentations maximales dans un  groupe  hermitien fournissaient les seuls exemples connus jusqu'ici où les invariants primitifs sont insuffisants. Cette note a  pour objet d'annoncer l'existence de nouvelles composantes  exotiques pour les espaces de modules pour les groupes $\SO(p,q)$, pour $2<p<q$.
\end{abstract}

\maketitle

\markleft{Aparicio-Arroyo, Bradlow, Collier, Garcia-Prada, Gothen, Oliveira}





\section{Introduction}

Representation varieties for closed oriented surfaces generally have more than one connected component. Some of the components are mundane in the sense that they are distinguished by obvious topological invariants and have no known special characteristics. Others are more alluring and unusual either because they are not detected by the primary invariants or because they have special geometric significance, or both. 

Instances of such `exotic' components are well understood in two situations. 
The first is for representations into a split real form of a complex semisimple group (see \cite{hitchin:1992}), in which case the exotic components are known as Hitchin components. 
The second occurs for representations into the isometry group for a non-compact Hermitian symmetric space (see \cite{bradlow-garcia-prada-gothen:2005}), in which case the components with so-called maximal Toledo invariant have exotic components. 
Both of these classes of exotic components have been called higher Teichm\"uller components in \cite{BIW2012} since they enjoy many of the geometric features of Teichm\"uller space.

The purpose of this note is to announce\footnotemark\footnotetext{Full
  details of the proof can be found in the preprint \cite{SOpq-complete}.} the existence of exotic components in the representation varieties for representations into $\SO(p,q)$, the special orthogonal groups with signature $(p,q)$. These groups are non-compact with two connected components, the component of the identity will be denoted by $\SO_0(p,q)$. The groups $\SO(p,q)$ are neither split nor of Hermitian type except in the special cases $p=2$ (in which case the groups are Hermitian) or the cases $q=p$ and $q=p+1$ (the split cases). The new $\SO(p,q)$ exotic components are thus not accounted for by previously known mechanisms.

For a closed surface $S$ and a Lie group $\G$, the representation variety $\Rep(S,\G)$ parameterizes the space of conjugacy classes of group homomorphisms $\rho:\pi_1(S)\to\G$. By the Non-Abelian Hodge (NAH) correspondence, the space $\Rep(S,\G)$ is homeomorphic to $\mathcal{M}(\Sigma,\G)$, the moduli space of polystable $\G$-Higgs bundles on a closed Riemann surface $\Sigma$ of the same genus as $S$. 
Our methods lie on the Higgs bundle side of the NAH correspondence; what we actually prove is the existence of exotic components for $\mathcal{M}(\Sigma,\SO(p,q))$. 

Though not accounted for by previously known mechanisms, evidence for these new exotic components has nevertheless steadily accumulated in recent years. 
The earliest indication came from Morse theoretic considerations based on the norm-squared of the Higgs field.  
This real-valued function defines a proper map and hence attains a local minimum on every connected component. 
The absolute minimum, i.e.\ the zero level, is attained on the components labeled by the primary topological invariants for $\SO(p,q)$-Higgs bundles. 
In her 2009 Ph.D. thesis \cite{aparicio-arroyo:2009} the first author described additional local minima at non-zero values.  
Since the moduli spaces are not smooth, Morse theoretic tools could not be employed directly to infer the existence of exotic components, but these exotic local minima revealed that possibility (now confirmed by the results in this note).

Further evidence came from the results in \cite{collier:2017.2} in the special case of $\SO(p,p+1)$.  Since these groups are split real forms, the moduli spaces have Hitchin components, but the results in \cite{collier:2017.2} show that these are not the only exotic components. 
For $p=2$, the fact that $\SO_0(2,3)$ is double covered by $\Sp(4,\R)$ leads to extra exotic components related to exotic components first detected in \cite{gothen} in the moduli spaces of $\Sp(4,\R)$-Higgs bundles. 
The results in \cite{collier:2017.2} show that these exotic components have counterparts for all $\SO(p,p+1)$.

Working on the other side of the Non-Abelian Hodge Correspondence, Guichard and Wienhard found reason to conjecture (see \cite{GW}) that the representation varieties $\Rep(S,\SO(p,q))$ should have special connected components not detected by the primary topological invariants. 
The conjectured components would parameterize representations characterized by a positivity condition which is a refinement of the Anosov condition introduced by Labourie \cite{AnosovFlowsLabourie}. 
The conjecture is based on the fact that apart from the split real forms and the real forms of Hermitian type, the only other non-exceptional groups which allow positive representations are the groups $\SO(p,q)$. 

The NAH correspondence is notoriously non-explicit, providing little guidance for matching individual Higgs bundles with specific surface group representations. 
In particular, except in special cases, it is difficult to identify Higgs bundles which correspond to representations satisfying the Guichard-Wienhard positivity condition. 
Nonetheless, with evidence coming from special cases, we conjecture that the Higgs bundles bundles in our exotic components all correspond to positive representations in the sense of Guichard-Wienhard, and hence that the exotic components we detect in $\mathcal{M}(\Sigma,\SO(p,q))$ correspond to the expected components in  $\Rep(S,\SO(p,q))$. In particular, the exotic components of $\Rep(S,\SO(p,q))$ would consist entirely of Anosov representations.

We note finally that very recently another facet of our exotic components has been revealed by Baraglia and Schaposnik in \cite{barraglia-schaposnik}. In this work, the space $\mathcal{M}(\Sigma,\SO(p,q))$ is regarded as a subvariety of $\mathcal{M}(\Sigma,\SO(p+q,\CC))$ and the intersection of $\mathcal{M}(\Sigma,\SO(p,q))$ with generic fibers of the Hitchin fibration is examined. 
Since their methods apply only to an open subset of $\mathcal{M}(\Sigma,\SO(p,q))$, they are not able to establish the existence of exotic components. They can however characterize the intersections of these components with generic fibers in terms of spectral data. 
The description has features in common with our description of the entire component but is not the same. It would be interesting to reconcile the two descriptions.

\section{The main result}

Let $S$ be a closed oriented surface of genus $g\geq 2$.  Fix a complex structure on $S$ and denote the resulting Riemann surface by $\Sigma$. Without loss of generality we assume that $p\leq q$.  While the concept of a $G$-Higgs bundle for non-compact real forms goes back to the pioneering work of Hitchin \cite{hitchin:1992}, the basic definitions and constructions have been elaborated on in several places including \cite{garcia-gothen-mundet:2008,schmitt:2008}. We give here only the essentials necessary to describe our results, and refer the reader to the cited references for more details.

\begin{definition} An $\SO(p,q)$-Higgs bundle on $\Sigma$ is defined by a triple $(V,W,\eta)$ where $V$ and $W$ are respectively rank $p$ and rank $q$ vector bundles with orthogonal structures such that\footnotemark\footnotetext{This condition is equivalent to $\det(W)\simeq\det(V)^{-1}$ since $\det(W)^2$ is trivial.} $\det(W)\simeq\det(V)$, and $\eta$ is a holomorphic bundle map $\eta:W\rightarrow V\otimes K$.  
\end{definition}

There are notions of stability, semistability, and polystability which apply to $\SO(p,q)$-Higgs bundles and which facilitate the construction of moduli spaces. We use the notation $\mathcal{M}(\SO(p,q))$ to denote the moduli space of polystable $\SO(p,q)$-Higgs bundles on $\Sigma$.

The cases with $p\leq2$ or $q\le p+1$ are somewhat special, so for clarity we state the main result without them, and discuss them in Section \ref{specialcases}. For $p>2$, principal $\mathrm{O}(p,\CC)$- bundles on $\Sigma$ are classified topologically by first and second Stiefel-Whitney classes, $sw_1\in H^1(S,\Z_2)$ and $sw_2\in H^2(S,\Z_2)$, coming from a reduction to the maximal compact subgroup $\mathrm{O}(p)$. These primary topological invariants are constant on connected components of the moduli space $\mathcal{M}(\SO(p,q))$.  Since $\det(W)\simeq\det(V)$, it follows that $sw_1(V)=sw_1(W)$.  The components of the moduli space $\mathcal{M}(\SO(p,q))$ are thus partially labeled by triples $(a,b,c)\in \Z_2^{2g}\times\Z_2\times\Z_2$, where 
\[\xymatrix@C=.3em{a=sw_1(V)\in H^1(\Sigma,\Z_2),&b=sw_2(V)\in H^2(\Sigma,\Z_2)&\text{and}&c=sw_2(W)\in H^2(\Sigma,\Z_2).}\]
Using the notation $\mathcal{M}^{a,b,c}(\SO(p,q))$ to denote the union of components labeled by $(a,b,c)$, we can thus write
\begin{equation}\label{Mpq-abc}
\mathcal{M}(\SO(p,q))=\coprod_{(a,b,c)\in \Z_2^{2g}\times\Z_2\times\Z_2}\mathcal{M}^{a,b,c}(\SO(p,q))\ .
\end{equation}






The stability notion for $\SO(p,q)$-Higgs bundles implies that a Higgs bundle $(V,W,\eta)$ with $\eta=0$ is polystable if and only if $V$ and $W$ are both polystable orthogonal bundles. This leads to the immediate identification of one connected component in each space $\mathcal{M}^{a,b,c}(\SO(p,q))$.  
We use the subscript `top' to designate these components, which contain $\SO(p,q)$-Higgs bundles with vanishing Higgs field.

\begin{proposition}\label{zerocomponents}  Assume that $2<p\le q$.  
For every $(a,b,c)\in \Z_2^{2g}\times\Z_2\times\Z_2$ the space $\mathcal{M}^{a,b,c}(\SO(p,q))$ has a non-empty connected component\footnotemark\footnotetext{In the case $p=2$ it is no longer true that $\mathcal{M}_{\rmin}^{a,b,c}(\SO(p,q))$ is non-empty for all $(a,b,c)$. In particular, if $a=0,$ then $V=L\oplus L^{-1}$ which (a) is polystable only if $\deg{L}=0$ and (b) has $sw_2(V)=\deg{L} \mod 2$. Thus $\mathcal{M}_{\rmin}^{0,b,c}(\SO(2,q))$ is empty if $b\ne 0$.}, denoted by $\mathcal{M}_{\rmin}^{a,b,c}(\SO(p,q))$, in which every point can be continuously deformed to the isomorphism class of an $\SO(p,q)$-Higgs bundle of the form $(V,W,\eta=0)$ where $V$ and $W$ are polystable orthogonal bundles.
\end{proposition}
We define
\begin{equation}
\mathcal{M}_{\rmin}(\SO(p,q))=\coprod_{a,b,c}  \mathcal{M}_{\rmin}^{a,b,c}(\SO(p,q))
\end{equation}
Our main result shows that the moduli space $\mathcal{M}(\SO(p,q))$ has additional `exotic' components disjoint from the components of $\mathcal{M}_{\rmin}(\SO(p,q))$.  We identify these exotic components as products of moduli spaces of so-called $L$-twisted Higgs bundles, where each factor $L$ is a positive power of the canonical bundle $K$.  

\begin{definition}  Let $L$ be a fixed holomorphic line bundle on $\Sigma$.  An $L$-twisted $\SO(1,n)$-Higgs bundle on $\Sigma$ is a triple $(I,W_0,\eta)$, where $W_0$ is a $\Or(n,\CC)$-bundle, $I$ is the rank one orthogonal bundle $det(W_0)$ and $\eta:W_0\to I\otimes L$ is a holomorphic bundle map.  
\end{definition}

\begin{remark} The bundle $I$ is a square root of the trivial bundle. If $I$ is trivial then the objects $(I,W_0,\eta)$ define $L$-twisted $\SO_0(1,n)$-Higgs bundles, where $\SO_0(1,n)$ denotes the connected component of the identity.
\end{remark}

\begin{remark} 
The notions of stability for Higgs bundles readily extend to $L$-twisted Higgs bundles and similarly allow the construction of moduli spaces. We use the notation $\mathcal{M}_L(\G)$ for the moduli space of polystable $L$-twisted $\G$-Higgs bundles. In particular, taking $L=K^p$, $\mathcal{M}_{K^p}(\SO(1,n))$ denotes the  moduli space of polystable $K^p$-twisted $\SO(1,n)$-Higgs bundles.  
\end{remark}

We get a decomposition similar to \eqref{Mpq-abc}, namely
\begin{equation}\label{M1n-ac}
\mathcal{M}_{K^p}(\SO(1,n))=\coprod_{(a,c)\in \Z_2^{2g}\times\Z_2}\mathcal{M}_{K^p}^{a,c}(\SO(1,n))\ ,
\end{equation}
\noi where $\mathcal{M}_{K^p}^{a,c}(\SO(1,n))$ denotes the component in which the $\SO(1,n)$-Higgs bundles are of the form $(I,W_0,\eta)$ with $a=sw_1(W_0)$ and $c=sw_2(W_0)$.

We can now state our main result.

\begin{theorem}[Main Theorem]\label{mainth} Fix integers $(p,q)$ such that $2< p< q-1$.  For each choice of $a\in \Z_2^{2g}$ and $c\in\Z_2$, the moduli space $\mathcal{M}(\SO(p,q))$ has a connected component disjoint from $\mathcal{M}_{\rmin}(\SO(p,q))$. This component is isomorphic to 
\begin{equation}\label{pqcayley}
\mathcal{M}_{K^p}^{a,c}(\SO(1,q-p+1))\times\mathcal{M}_{K^2}(\SO_0(1,1))\times\cdots\times\mathcal{M}_{K^{2p-2}}(\SO_0(1,1))~,
\end{equation}
and lies in the sector $\mathcal{M}^{\alpha,0,c}(\SO(p,q))$ where $\alpha=a$ if $p$ is odd and $\alpha=0$ if $p$ is even. Moreover, $\cM(\SO(p,q))$ has no other connected components. 
 \end{theorem}

\begin{remark}
The group $\SO_0(1,1)$ is the connected component of the identity in $\SO(1,1)$, and $\mathcal{M}_{K^{2j}}(\SO_0(1,1))$ can be identified with $H^0(K^{2j})$. Thus, we can replace \eqref{pqcayley} with
\begin{equation}\label{pqcayley.2}
\mathcal{M}_{K^p}^{a,c}(\SO(1,q-p+1))\times \bigoplus_{j=1}^{p-1}H^0(K^{2j})~.
\end{equation}
\end{remark}

\begin{remark}
The existence of the exotic components described by \eqref{pqcayley} was proven for $p=2$ in \cite{bradlow-garcia-prada-gothen:2005}. They are the exotic components with maximal Toledo invariant arising from Cayley correspondence (see Section \ref{SO2q}). 
In particular, Theorem \ref{mainth} can be viewed as a generalized Cayley correspondence. Contrary to the cases $p>2$, there are components of $\mathcal{M}(\SO(2,q))$ which are not in the family described by the theorem and also not in $\mathcal{M}_{\rmin}(\SO(2,q))$. These are the components with non-maximal and non-zero Toledo invariant. 
\end{remark}



For $2<p<q-1$ and each $(a,c)$, we show that the space $\cM^{a,c}_{K^p}(\SO(1,q-p+1))$ is connected. As an immediate corollary, this gives a count of the connected components of $\cM(\SO(p,q)).$
\begin{corollary} For $2<p<q-1$, the moduli space $\mathcal{M}(\SO(p,q))$ has $3\times 2^{2g+1}$ connected components, $2^{2g+1}$ of which are exotic components disjoint from $\mathcal{M}_\mathrm{top}(\SO(p,q))$. 
\end{corollary}

\begin{remark} As a further consequence of Theorem \ref{mainth} we have the following dichotomy for a Higgs bundle in $\mathcal{M}(\SO(p,q))$: either the Higgs field can be deformed to zero, or the Higgs bundle is determined by a polystable $\mathrm{O}(q-p+1,\CC)$-bundle and a Higgs bundle in a Hitchin component for $\SO(p,p-1)$.
\end{remark}





\section{Special cases}\label{specialcases}

We now describe the special cases not covered by Theorem \ref{mainth}. Even though they are mainly accounted for by previously known phenomena, we will see in the next section that all cases fit within a unified framework.
\subsection{The case $q=p+1$}  If $q=p+1$, then $\mathcal{M}_{K^p}^{a,c}(\SO(1,q-p+1))=\mathcal{M}_{K^p}^{a,c}(\SO(1,2))$, which is not always connected.  Indeed, if $a=0$, then the Higgs bundles represented in $\mathcal{M}_{K^p}^{0,c}(\SO(1,2))$ can be taken to be of the form $(\mathcal{O}, L\oplus L^{-1}, \eta)$, where $L$ is a non-negative degree $d$ line bundle. Stability considerations impose a bound on $d$ so that 
\begin{equation}\label{so120}
\mathcal{M}_{K^p}^{0,c}(\SO(1,2))=\coprod_{\substack{0\le d\le p(2g-2)\\ d=c\,(\mathrm{mod}\ 2)}}\mathcal{M}_{K^p}^{d}(\SO(1,2)).
\end{equation}
Moreover, (see \cite[Theorem 1]{collier:2017.2}) for each integer $d\in (0, 2g-2]$, the moduli space $\mathcal{M}_{K^p}^{d}(\SO(1,2))$ is diffeomorphic to a vector bundle of rank $d+g-1$ over the $(2g-2-d)^{th}$-symmetric product $\Sym^{2g-2+d}(\Sigma)$. In particular, the components $\mathcal{M}_{K^p}^{d}(\SO(1,2))$ are smooth and connected.  

The moduli spaces $\mathcal{M}(\SO(p,p+1))$ have been analyzed in \cite{collier:2017.2}. It was shown there that the topological invariants for $\SO(p,p+1)$-Higgs bundles, i.e.\  the triples $(a,b,c)$, do not distinguish all connected components.  Two families of exotic components were identified. The components in the first family are labeled by an integer, $d$, in the range $0\le d\le p(2g-2)$, while those in the second family are labeled by a pair $(a,c)\in(\Z^{2g}-\{0\})\times\Z_2$. Though not described in this way in \cite{collier:2017.2}, these families can be identified as follows:

\begin{itemize}
\item In the the family labeled by $d$, each member is isomorphic to
\begin{equation}
\mathcal{M}_{K^p}^{d}(\SO(1,2))\times\mathcal{M}_{K^2}(\SO_0(1,1))\times\cdots\times\mathcal{M}_{K^{2p-2}}(\SO_0(1,1)),
\end{equation}
\noi where $\mathcal{M}_{K^p}^{d}(\SO(1,2))$ is one of the components of $\mathcal{M}_{K^p}^{0,c}(\SO(1,2))$ as in \eqref{so120}.

\item In the family labeled by $(a,c)$, each member is isomorphic to
\begin{equation}
\mathcal{M}_{K^p}^{a,c}(\SO(1,2))\times\mathcal{M}_{K^2}(\SO_0(1,1))\times\cdots\times \mathcal{M}_{K^{2p-2}}(\SO_0(1,1)).
\end{equation}
\end{itemize}

\noi The components are thus precisely those identified by Theorem \ref{mainth} in the case  $q=p+1$.  The component count in this case is, however, different from the case $q>p+1$.

\begin{corollary} For $p>2$, the moduli space $\mathcal{M}(\SO(p,p+1))$ has $3\times 2^{2g+1}+2p(g-1)-1$ connected components.  Among those, there are $2^{2g+1}+2p(g-1)-1$ `exotic' components which are disjoint from $\cM_{\rmin}(\SO(p,p+1))$.  
\end{corollary}

\subsection{The case $q=p$} In this case $\mathcal{M}_{K^p}(\SO(1,q-p+1))=\mathcal{M}_{K^p}(\SO(1,1))$.  A $K^p$-twisted $\SO(1,1)$-Higgs bundle consists of a triple $(I,I,\eta)$ where $I$ is a square root of the trivial bundle $\mathcal{O}$ and $\eta\in H^0(K^p)$.  Such Higgs bundles are labeled by a single Stiefel-Whitney class, namely $a=sw_1(I)$, so that

\begin{equation}
\mathcal{M}_{K^p}(\SO(1,1))=\coprod_{a\in H^1(\Sigma,\Z_2)}\mathcal{M}^a_{K^p}(\SO(1,1)).
\end{equation}

\noi With $q=p,$ Theorem \ref{mainth} thus gives $2^{2g}$ exotic components of $\mathcal{M}(\SO(p,p))$ isomorphic to the moduli spaces
\begin{equation}
\mathcal{M}_{K^p}^{a}(\SO(1,1))\times\mathcal{M}_{K^2}(\SO_0(1,1))\times\cdots\times\mathcal{M}_{K^{2p-2}}(\SO_0(1,1)).
\end{equation}

\noi For each $a$, we can identify $\mathcal{M}_{K^p}^{a}(\SO(1,1))$ with $H^0(K^p)$. Thus, each exotic component is isomorphic to $H^0(K^p)\oplus\bigoplus\limits_{j=1}^{p-1}(H^0(K^{2j})$.  This recovers 
 the Hitchin component in $\mathcal{M}(\SO_0(p,p))$ when $a=0$.

\subsection{The case $p=2<q$.}\label{SO2q} An $\SO(2,q)$-Higgs bundle is defined by a triple $(V,W,\eta)$ in which $V$ is an $\mathrm{O}(2,\CC)$-bundle.  If $sw_1(V)=0$, i.e.\ if the structure group of $V$ reduces to $\SO(2,\CC)$, then $V$ can be assumed to be a direct sum of line bundles of the form $V=L\oplus L^{-1}$, with orthogonal structure given $q_V=\smtrx{0&1\\1&0}$ in this splitting. Note that the second Stiefel-Whitney class of the orthogonal bundle $L\oplus L^{-1}$ is given by $sw_2=d \,(\mathrm{mod}\ 2)$ where $d=\deg(L)\geq0$.  

For the groups $\SO(2,q)$, the connected components of the identity are isometry groups of Hermitian symmetric spaces of non-compact type. As such, the Higgs bundles  have an associated Toledo invariant which, up to a normalization constant, is integer-valued but subject to a Milnor-Wood bound (see \cite{BGR,bradlow-garcia-prada-gothen:2005}). For an $\SO_0(2,q)$-Higgs bundle $(L\oplus L^{-1}, W, \eta)$, the Toledo invariant can be identified as the degree $d$ of $L$ and the Milnor-Wood bound\footnote{For the group $\SO_0(2,q)$ the Milnor-Wood inequality is $2-2g\leq d\leq 2g-2$. However, the automorphism switching the sign of the Toledo invariant (which is an outer automorphism for $\SO_0(2,q)$) can be realized as an $\SO(2,q)$ inner automorphism.} is $0\leq d\le 2g-2$.
We thus get
\begin{equation}
\mathcal{M}^{0,b,c}(\SO(2,q))=\coprod_{\substack{0\leq d\leq 2g-2\\ d=b \,(\mathrm{mod}\ 2)}}\mathcal{M}^{d,c}(\SO_0(2,q)),
\end{equation}
\noi where $\mathcal{M}^{d,c}(\SO(2,q))$ denotes the component in which $\deg(L)=d$ and $sw_2(W)=c$.  

The components where $d=2g-2$ specializes further because in these components:
\begin{enumerate}
\item $L$ has to be isomorphic to $KI$ where $I^2=\mathcal{O}$, and
\item $W$ decomposes as $W=I\oplus W_0$ where $W_0$ is a rank $q-1$ orthogonal bundle with $sw_1(W_0)=I$.
\end{enumerate}

\noi  As shown in \cite{CTT}, an $\SO(2,q)$-Higgs bundle with $L=KI$,  $W=I\oplus W_0$ and $\eta=[q_2,\beta]:I\oplus W_0\rightarrow KI$ is defined by a $K^2$-twisted $\SO(1,q-1)$-Higgs bundle $(I,W_0,\beta)$ together with a quadratic differential $q_2$. Denoting $\mathcal{M}^{2g-2,c}(\SO(2,q))$ by $\mathcal{M}^{\mathrm{max},c}(\SO(2,q))$ it follows  that 
\begin{equation}
\mathcal{M}^{\mathrm{max},c}(\SO(2,q))=\coprod_{a\in H^1(\Sigma,\Z_2)}\mathcal{M}^{a,c}_{K^2}(\SO(1,q-1))\times H^0(K^2)=\coprod_{a}\mathcal{M}^{a,c}_{K^2}(\SO(1,q-1)\times \SO_0(1,1))\end{equation}
\noi where $a=sw_1(I)$.  We thus get
\begin{equation}
\mathcal{M}^{0,0,c}(\SO(2,q))=\coprod_{\substack{0\leq d<2g-2\\ d=0 \,(\mathrm{mod}\ 2)}}\mathcal{M}^{d,c}(\SO(2,q))\ \sqcup\ \coprod_{a}\mathcal{M}^{a,c}_{K^2}(\SO(1,q-1)\times \SO_0(1,1)).
\end{equation}

\noi The group $\SO(1,q-1)\times\SO_0(1,1)$ is known as the Cayley partner to $\SO(2,q)$ and the objects in $\mathcal{M}^{a,c}_{K^2}(\SO(1,q-1)\times \SO_0(1,1))$ are called Cayley partners to the Higgs bundles in $\mathcal{M}^{\mathrm{max},c}(\SO(2,q))$.  Such Cayley partners are known to emerge in maximal components of $\mathcal{M}(\G)$ whenever $\G$ is the isometry group of a Hermitian symmetric space of tube type (see \cite{BGR,bradlow-garcia-prada-gothen:2005}). Comparing to Theorem \ref{mainth}, we see that the exotic components in $\mathcal{M}(\SO(p,q))$ are direct generalizations of these Cayley partners to the maximal components $\mathcal{M}^{\mathrm{max},c}(\SO(2,q))$.

%
%
%
%
%
%
%

\section{What we actually prove -- the main ideas}

Theorem \ref{mainth} shows not only that additional exotic components exist, but also gives a model which describes them. Indeed, given the model, the result is proved directly by constructing a suitable map from the model to the moduli space $\mathcal{M}(\SO(p,q))$. The model is itself built from moduli spaces (of $K^j$-twisted Higgs bundles), so that in both the domain and target of our map the points represent equivalence classes of objects. We first describe a map between the objects, and then show that it descends to the appropriate moduli spaces where it defines a homeomorphism onto a connected component.

The map relies in part on a parameterization of the Hitchin components of the moduli spaces $\mathcal{M}(\SO(p-1,p))$.  Viewing theses moduli spaces as subspaces of $\mathcal{M}(\SO(2p-1,\C))$, the parameterization is given by a section of the Hitchin fibration for $\mathcal{M}(\SO(2p-1,\C))$.  The fibration is defined by $\SO(2p-1,\C)$-invariant polynomials evaluated on the Higgs field, giving a map to $\bigoplus\limits_{j=1}^{p-1}H^0(K^{2j})$. Hitchin showed in \cite{hitchin:1992} that the fibration admits sections which parameterizes connected components of $\mathcal{M}(\SO(p-1,p))\subset\mathcal{M}(\SO(2p-1),\C)$. 

The $\SO(p,p-1)$-Higgs bundles in the image of the section can be taken to be of the form $(\mathcal{K}_{p},\mathcal{K}_{p-1},\sigma(q_2,q_4,\dots,q_{2p-2}))$ where
\begin{equation}\label{Eq cal K notation}
\mathcal{K}_{p}=K^{p-1}\oplus K^{p-3}\oplus\cdots\oplus K^{1-p},
\end{equation}
and $\sigma$ is given by a map
\begin{equation}
	\label{EQ Hitchin section map}
	\sigma:\bigoplus_{j=1}^nH^0(K^{2j})\rightarrow\Hom(\mathcal{K}_{p},\mathcal{K}_{p-1})\otimes K.
\end{equation}


\noi Our Main Theorem is a consequence of the following:

\begin{theorem}\label{mainmetatheorem}
 Let $(I,W_0,\eta_p)$ be a $K^p$-twisted $\SO(1,q-p+1)$-Higgs bundle and take differentials $q_{2j}\in  H^0(K^{2j})$ for $j=1,\dots p-1$. Using the notation from \eqref{Eq cal K notation}, consider the $\SO(p,q)$-Higgs bundle 
$(V,W,\eta)$ defined by
\begin{equation}
	\label{EQ etaform}\xymatrix{V=I\otimes \mathcal{K}_p,& W=W_0\oplus I\otimes\mathcal{K}_{p-1}&\text{and}}\ \ \ \ \eta=\mtrx{\mu&\sigma(\ \vec{q}\ )}:W\to V\otimes K,
\end{equation}
where $\mu=\smtrx{\eta_p\\0\\\vdots\\0}:W_0\to \mathcal{K}_p\otimes I\otimes K,$ $\vec{q}=(q_2,\dots,q_{2p-2})$ and $\sigma$ is the map from \eqref{EQ Hitchin section map}.
\begin{enumerate}
\item $(V,W,\eta)$ is polystable if and only if 
 $(I,W_0,\mu)$ is  polystable. The map
\[
 ((I,W_0,\mu); q_2,q_4,\dots,q_{2p-2})\longmapsto (V,W,\eta)\] thus descends to define a map
\[\Psi:\xymatrix{\cM_{K^p}(\SO(1,q-p+1))\times\displaystyle\prod\limits_{i=1}^{p-1} \cM_{K^{2i}}(\SO_0(1,1))\ar[rr]&& \cM(\SO(p,q))}.\]

\item The map $\Psi$ is injective.

\item The map $\Psi$ is open and closed.
\item The image of $\Psi$ is disjoint from the components $\mathcal{M}_0(\SO(p,q))$.

\end{enumerate}
\end{theorem}

After showing that the $\SO(p,q)$-Higgs bundle \eqref{EQ etaform} is polystable if and only if $(I,W_0,\eta_p)$ is a polystable $K^p$-twisted $\SO(1,q-p+1)$-Higgs bundle, it is shown that two Higgs bundles  of the form \eqref{EQ etaform} lie in the same gauge orbit if and only $(I,W_0,\eta_p)$ lie in the same gauge orbit. This proves the map $\Psi$ on moduli spaces is injective. 

To prove that the map $\Psi$ is open, we analyze the local structures of both the domain and target of the map $\Psi.$ In particular, we show that the second hyper-cohomology groups in the deformation complexes vanish at all points of both the domain and image of the map. This technical result together with an appropriately equivariant isomorphism of the first hyper-cohomologies of the deformation complexes leads to openness of the map $\Psi.$






For the closedness of the map $\Psi$, we prove the contrapositive, i.e.\ we show that if $\Psi$ fails to be a closed mapping, then we can find a divergent sequence of points $\{x_i\}$ in the domain 
such that $\{\Psi(x_i)\}$ converges in $\mathcal{M}(\SO(p,q))$.  We use the properness of the Hitchin fibrations for $\cM_{K^p}(\SO(1,q-p+1))$ and $\mathcal{M}(\SO(p,q))$ to show that no such sequence exists.  In particular, since the fibers of the Hitchin fibration are compact, if the sequence $\{x_i\}$ diverges, then the projection of the sequence onto the base of the fibration, say $\{y_i\}$, must diverge.  Using the map induced on the Hitchin bases by $\Psi$, say $\hat{\Psi}$, if $\{y_i\}$ diverges, then so does the sequence $\{\hat{\Psi}(y_i)\}$, and hence so does the sequence $\{\Psi(x_i)\}$. 

Parts (1)-(3) prove that the image of $\Psi$ forms a union of connected components. We show further that any Higgs bundle not in the image of $\Psi$ can be deformed to one in which the Higgs field is zero. Since, by construction,  the Higgs field never vanishes in the image of $\Psi$, this proves (4).

\section{Conjecture}
Representations in the so-called higher Teichm\"uller components of $\Rep(S,\G)$ are examples of the important class of Anosov representations introduced by Labourie in \cite{AnosovFlowsLabourie}. As a result, they have many interesting geometric and dynamical properties, One important property of the set of Anosov representations is that they are open in the representation variety. They are however not closed in general. 

Recently, Guichard and Wienhard \cite{GW} introduced the notion of a positive Anosov representation which refines the notion of an Anosov representation. Here positivity depends on a choice of parabolic subgroup of $\G,$ and like Anosov representations, positive representations are open in the representation variety. Guichard and Wienhard conjecture that the set of positive representations is also closed in the representation variety. As a result, positive representations define connected components of the representation variety which consist entirely of Anosov representations. 

For the two known families of higher Teichm\"uller components, this has been established. Namely, Hitchin representations are positive with respect to the Borel subgroup \cite{AnosovFlowsLabourie,fock_goncharov_2006} and, for a Hermitian group $\G$ of tube type, maximal representations are positive with respect to the parabolic subgroup which gives rise to the Shilov boundary of the Riemannian symmetric space of $\G$ \cite{MaxRepsAnosov}.
Interestingly, for $p>q$, $\SO(p,q)$ also admits a notion of positivity with respect to the generalized flag variety consisting of flags 
$V_1\subset\cdots\subset V_{p- 1} \subset \R^{p+q}$ where $V_j$ is an isotropic (with respect to a signature $(p,q)$ inner product) $j$-plane.  

It is natural to conjecture that the exotic connected components of Theorem \ref{mainth} correspond to connected components of positive representations in $\SO(p,q).$
\begin{conjecture}Under the NAH correspondence, the exotic components identified in Theorem \ref{mainth} correspond to the positive components conjectured to exist by Guichard-Wienhard.
\end{conjecture}

It is shown in \cite{collier:2017.2} that many of the exotic components of $\mathcal{M}(\SO(p,p+1))$ contain points corresponding to positive representations. Generalizing these techniques, we can show that for $q>p+1$ every exotic components of $\mathcal{M}(\SO(p,q))$ contains Higgs bundles which correspond to positive representations.  Thus, if Guichard and Wienhard's conjecture on closedness of positive representations is true, it would imply the above conjecture and Theorem \ref{mainth} would give a count of the connected components of positive representations.

\section{Acknowledgements}
We are pleased to acknowledge many enlightening conversations over several years with Anna Wienhard and Olivier Guichard  about positivity and its likely implications for $\SO(p,q)$-moduli spaces.  We also thank the referee for a careful reading of this note and helpful comments.

The authors acknow\-ledge support from U.S. National Science Foundation grants DMS 1107452, 1107263, 1107367 ``RNMS: GEometric structures And Representation varieties" (the GEAR Network). The fourth author was partially supported by the Spanish MINECO under ICMAT Severo Ochoa project No. SEV-2015-0554, and under grant No. MTM2013-43963-P. The fifth and sixth authors were partially supported by CMUP (UID/MAT/00144/2013) and the project PTDC/MAT-GEO/2823/2014 funded by FCT (Portugal) with national funds. The sixth author was also partially supported by the Post-Doctoral fellowship SFRH/BPD/100996/2014 funded by FCT (Portugal) with national funds.

\end{document}